\documentclass[a4paper,oneside,10pt]{article}%
\usepackage{amsmath}
\usepackage{amsfonts}
\usepackage{amssymb}
\usepackage{graphicx}
\usepackage{color}
\usepackage[square,numbers,sort&compress]{natbib}%
\providecommand{\U}[1]{\protect \rule{.1in}{.1in}}

\newtheorem{theorem}{Theorem}[section]

\newtheorem{definition}[theorem]{Definition}

\newtheorem{lemma}[theorem]{Lemma}

\newtheorem{proposition}[theorem]{Proposition}
\newtheorem{remark}[theorem]{Remark}

\newenvironment{proof}[1][Proof]{\noindent \textbf{#1.} }{\  $\Box$}
\numberwithin{equation}{section}

\begin{document}

\title{Dynamic programming principle and Hamilton-Jacobi-Bellman equation under
nonlinear expectation}
\author{Mingshang Hu \thanks{Zhongtai Securities Institute for Financial Studies,
Shandong University, Jinan, Shandong 250100, PR China. humingshang@sdu.edu.cn.
Research supported by National Key R\&D Program of China (No. 2018YFA0703900)
and NSF (No. 11671231). }
\and Shaolin Ji\thanks{Zhongtai Securities Institute for Financial Studies,
Shandong University, Jinan, Shandong 250100, PR China. jsl@sdu.edu.cn.
Research supported by NSF (No. 11971263 and 11871458). }
\and Xiaojuan Li\thanks{Zhongtai Securities Institute for Financial Studies,
Shandong University, Jinan 250100, China. Email: lixiaojuan@mail.sdu.edu.cn.} }
\maketitle

\textbf{Abstract}. In this paper, we study a stochastic recursive optimal
control problem in which the value functional is defined by the solution of a
backward stochastic differential equation (BSDE) under $\tilde{G}%
$-expectation. Under standard assumptions, we establish the comparison theorem
for this kind of BSDE and give a novel and simple method to obtain the dynamic
programming principle. Finally, we prove that the value function is the unique
viscosity solution of a type of fully nonlinear HJB equation.



\addcontentsline{toc}{section}{\hspace*{1.8em}Abstract}

\section{Introduction}

Motivated by the model uncertainty in finance, Peng \cite{Peng2005, P07a,
P08a} established the theory of $G$-expectation which is a consistent
sublinear expectation and does not require a probability space. The
representation of $G$-expectation as the supremum of expectations over a set
of nondominated probability measures was obtained in \cite{DHP11, HP09}. Due
to this set of nondominated probability measures, the backward stochastic
differential equation (BSDE for short) is completely different from the
classical one. Hu et al. \cite{HJPS1} obtained an existence and uniqueness
theorem for a new kind of BSDE driven by $G$-Brownian motion. In addition,
there are other advances in this direction. Denis and Martini
\cite{DenisMartini2006} developed quasi-sure stochastic analysis. Soner et al.
\cite{STZ11} obtained an existence and uniqueness theorem for a new type of
BSDE ($2$BSDE) under a family of nondominated probability measures.

Recently, Hu and Ji \cite{HJ1} studied the following stochastic recursive
optimal control problem under $G$-expectation:%
\begin{equation}
\left \{
\begin{array}
[c]{rl}%
dX_{s}^{t,x,u}= & b(s,X_{s}^{t,x,u},u_{s})ds+h_{ij}(s,X_{s}^{t,x,u}%
,u_{s})d\langle B^{i},B^{j}\rangle_{s}+\sigma(s,X_{s}^{t,x,u},u_{s})dB_{s},\\
X_{t}^{t,x,u}= & x,
\end{array}
\right.  \label{ine1}%
\end{equation}%
\begin{equation}%
\begin{array}
[c]{rl}%
Y_{s}^{t,x,u}= & \Phi(X_{T}^{t,x,u})+\int_{s}^{T}f(r,X_{r}^{t,x,u}%
,Y_{r}^{t,x,u},Z_{r}^{t,x,u},u_{r})dr+\int_{s}^{T}g_{ij}(r,X_{r}^{t,x,u}%
,Y_{r}^{t,x,u},Z_{r}^{t,x,u},u_{r})d\langle B^{i},B^{j}\rangle_{r}\\
& -\int_{s}^{T}Z_{r}^{t,x,u}dB_{r}-(K_{T}^{t,x,u}-K_{s}^{t,xu}),\text{ }%
s\in \lbrack t,T].
\end{array}
\label{ine2}%
\end{equation}
The value function is defined as%
\begin{equation}
V(t,x):=\underset{u\in \mathcal{U}[t,T]}{ess\inf}Y_{t}^{t,x,u}. \label{ine3}%
\end{equation}

As pointed out in \cite{HJ1}, the value function defined in (\ref{ine3}) is a
$\inf \sup$ problem, which is known as the robust optimal control problem. For
recent development of robust control problem under a set of nondominated
probability measures, we refer the readers to \cite{DK, EJ-1, EJ-2, MPZ, PZ}
and the references therein. When $G$ is linear, the above optimal control
problem is classical stochastic recursive optimal control problem, which was
first studied by Peng in \cite{peng-dpp}. For the development of classical
stochastic recursive optimal control problem, we refer the readers to
\cite{BH, BL, EPQ, HJX, LW, MPY, MY, Tang, WY, J.Yong} and the references therein.

The nonlinear part with respect to $\partial_{xx}^{2}V$ in the HJB equation
related to the optimal control problem (\ref{ine1}) and (\ref{ine2}) is the
$\inf \sup$ of a family of linear part with respect to $\partial_{xx}^{2}V$. Up
to our knowledge, this $\inf \sup$ representation is the only result that has
been made so far in the optimal control problem. In order to obtain the fully
nonlinear representation, we want to study the stochastic recursive optimal
control problem under $\tilde{G}$-expectation. Here $\tilde{G}$ is any
function dominated by $G$ in the meaning of (\ref{newe1}). More precisely, we
consider the following BSDE under $\tilde{G}$-expectation:%
\begin{equation}
Y_{s}^{t,x,u}=\mathbb{\tilde{E}}_{s}\left[  \Phi(X_{T}^{t,x,u})+\int_{s}%
^{T}f(r,X_{r}^{t,x,u},Y_{r}^{t,x,u},u_{r})dr+\int_{s}^{T}g_{ij}(r,X_{r}%
^{t,x,u},Y_{r}^{t,x,u},u_{r})d\langle B^{i},B^{j}\rangle_{r}\right]  .
\label{ine4}%
\end{equation}
The new optimal control problem is (\ref{ine1}) and (\ref{ine4}), and the
value function is still defined as (\ref{ine3}). It is worth pointing out that
the BSDE (\ref{ine4}) under $\tilde{G}$-expectation does not contain $Z$,
which is an important open problem.

In this paper, we study the dynamic programming principle (DPP) and HJB
equation for optimal control problem (\ref{ine1}) and (\ref{ine4}). Firstly,
we establish the comparison theorem for BSDE (\ref{ine4}), which is new in the
literature. Secondly, for each $\xi \in L_{G}^{2}(\Omega)$, we prove that there
exists a sequence of simple random variables $\xi_{k}\in L_{G}^{2}(\Omega)$
such that $\xi_{k}$ converges to $\xi$ in the sense of $L_{G}^{2}$. Based on
this approximation, we give a new method to prove the DPP, which still holds
for the optimal control problem (\ref{ine1}) and (\ref{ine2}) and is easier
than the implied partition method in \cite{HJ1}. At last, we prove that $V$ is
the unique viscosity solution of a type of fully nonlinear HJB equation, which
is not the $\inf \sup$ representation with respect to $\partial_{xx}^{2}V$.

This paper is organized as follows. We recall some basic results on
$G$-expectation and $\tilde{G}$-expectation in Section 2. In Section 3, we
formulate our stochastic recursive optimal control problem under $\tilde{G}%
$-expectation. In Section 4, we prove the properties of the value function and
obtain the DPP. We prove that the value function is the unique viscosity
solution of a type of fully nonlinear HJB equation in Section 5.

\section{Preliminaries}

Let $T>0$ be given and let $\Omega_{T}=C_{0}([0,T];\mathbb{R}^{d})$ be the
space of $\mathbb{R}^{d}$-valued continuous functions on $[0,T]$ with
$\omega_{0}=0$. The canonical process $B_{t}(\omega):=\omega_{t}$, for
$\omega \in \Omega_{T}$ and $t\in \lbrack0,T]$. Set%
\[
Lip(\Omega_{T}):=\{ \varphi(B_{t_{1}},B_{t_{2}}-B_{t_{1}},\ldots,B_{t_{N}%
}-B_{t_{N-1}}):N\geq1,t_{1}<\cdots<t_{N}\leq T,\varphi \in C_{b.Lip}%
(\mathbb{R}^{d\times N})\},
\]
where $C_{b.Lip}(\mathbb{R}^{d\times N})$ denotes the space of bounded
Lipschitz functions on $\mathbb{R}^{d\times N}$.

Let $G:\mathbb{S}_{d}\rightarrow \mathbb{R}$ be a given monotonic and sublinear
function, where $\mathbb{S}_{d}$ denotes the set of $d\times d$ symmetric
matrices. Peng \cite{P07a, P08a} constructed a $G$-expectation space
$(\Omega_{T},Lip(\Omega_{T}),\mathbb{\hat{E}},(\mathbb{\hat{E}}_{t}%
)_{t\in \lbrack0,T]})$, which is a consistent sublinear expectation space. The
canonical process $(B_{t})_{t\in \lbrack0,T]}$ is called $G$-Brownian motion
under $\mathbb{\hat{E}}$. Throughout this paper, we suppose that $G$ is
non-degenerate, i.e., there exists a $\underline{\sigma}^{2}>0$ such that
$G(A)-G(B)\geq \frac{1}{2}\underline{\sigma}^{2}\mathrm{tr}[A-B]$ for any
$A\geq B$. Furthermore, let $\tilde{G}:\mathbb{S}_{d}\rightarrow \mathbb{R}$ be
any given monotonic function dominated by $G$, i.e., for $A_{1}$, $A_{2}%
\in \mathbb{S}_{d}$,%

\begin{equation}
\left \{
\begin{array}
[c]{l}%
\tilde{G}(0)=0,\\
\tilde{G}(A_{1})\geq \tilde{G}(A_{2})\text{ if }A_{1}\geq A_{2},\\
\tilde{G}(A_{1})-\tilde{G}(A_{2})\leq G(A_{1}-A_{2}).
\end{array}
\right.  \label{newe1}%
\end{equation}
Peng also constructed a $\tilde{G}$-expectation space $(\Omega_{T}%
,Lip(\Omega_{T}),\mathbb{\tilde{E}},(\mathbb{\tilde{E}}_{t})_{t\in \lbrack
0,T]})$ in \cite{Peng2005, P2019}, which is a consistent nonlinear expectation
space satisfying%
\begin{equation}
\mathbb{\tilde{E}}_{t}[X]-\mathbb{\tilde{E}}_{t}[Y]\leq \mathbb{\hat{E}}%
_{t}[X-Y]\text{ for }X,Y\in Lip(\Omega_{T}),\text{ }t\in \lbrack0,T].
\label{e1}%
\end{equation}

Denote by $L_{G}^{p}(\Omega_{T})$ the completion of $Lip(\Omega_{T})$ under
the norm $||X||_{L_{G}^{p}}:=(\mathbb{\hat{E}}[|X|^{p}])^{1/p}$ for $p\geq1$.
For each $t\in \lbrack0,T]$, the conditional $G$-expectation and $\tilde{G}%
$-expectation can be continuously extended to $L_{G}^{1}(\Omega_{T})$ under
the norm $||\cdot||_{L_{G}^{1}}$, and still satisfy the relation (\ref{e1})
for $X,Y\in L_{G}^{1}(\Omega_{T})$.

\begin{definition}
Let $M_{G}^{0}(0,T)$ be the space of simple processes in the following form:
for each $N\in \mathbb{N}$ and $0=t_{0}<\cdots<t_{N}=T$,%
\[
\eta_{t}=\sum_{k=0}^{N-1}\xi_{k}I_{[t_{k},t_{k+1})}(t),
\]
where $\xi_{k}\in Lip(\Omega_{t_{k}})$ for $k=0,1,\ldots,N-1$.
\end{definition}

Denote by $M_{G}^{p}(0,T)$ the completion of $M_{G}^{0}(0,T)$ under the norm
$||\eta||_{M_{G}^{p}}:=(\mathbb{\hat{E}}[\int_{0}^{T}|\eta_{t}|^{p}dt])^{1/p}$
for $p\geq1$. For each $\eta^{k}\in M_{G}^{2}(0,T)$, $k=1,\ldots,d$, denote
$\eta=(\eta^{1},\ldots,\eta^{d})^{T}\in M_{G}^{2}(0,T;\mathbb{R}^{d})$, the
$G$-It\^{o} integral $\int_{0}^{T}\eta_{t}^{T}dB_{t}$ is well defined, see
Peng \cite{P07a, P08a, P2019}.

\begin{theorem}
(\cite{DHP11, HP09}) There exists a weakly compact set of probability measures
$\mathcal{P}$ on $(\Omega_{T},\mathcal{B}(\Omega_{T}))$ such that%
\[
\mathbb{\hat{E}}[X]=\sup_{P\in \mathcal{P}}E_{P}[X]\text{ for all }X\in
L_{G}^{1}(\Omega_{T}).
\]
$\mathcal{P}$ is called a set that represents $\mathbb{\hat{E}}$.
\end{theorem}

For this $\mathcal{P}$, we define capacity%
\[
c(A):=\sup_{P\in \mathcal{P}}P(A)\text{ for }A\in \mathcal{B}(\Omega_{T}).
\]
A set $A\in \mathcal{B}(\Omega_{T})$ is polar if $c(A)=0$. A property holds
\textquotedblleft quasi-surely" (q.s. for short) if it holds outside a polar
set. In the following, we do not distinguish two random variables $X$ and $Y$
if $X=Y$ q.s.

\section{Stochastic optimal control problem}

Let $U$ be a given compact set of $\mathbb{R}^{m}$. For each $t\in \lbrack
0,T]$, we denote by
\[
\mathcal{U}[t,T]:=\{u:u\in M_{G}^{2}(t,T;\mathbb{R}^{m})\text{ with values in
}U\}
\]
the set of admissible controls on $[t,T]$.

In the following, we use Einstein summation convention. For each given
$t\in \lbrack0,T]$, $\xi \in L_{G}^{p}(\Omega_{t};\mathbb{R}^{n})$ with $p\geq2$
and $u\in \mathcal{U}[t,T]$, we consider the following forward and backward
SDEs:%
\begin{equation}
\left \{
\begin{array}
[c]{rl}%
dX_{s}^{t,\xi,u}= & b(s,X_{s}^{t,\xi,u},u_{s})ds+h_{ij}(s,X_{s}^{t,\xi
,u},u_{s})d\langle B^{i},B^{j}\rangle_{s}+\sigma(s,X_{s}^{t,\xi,u}%
,u_{s})dB_{s},\\
X_{t}^{t,\xi,u}= & \xi,
\end{array}
\right.  \label{e2}%
\end{equation}
and%
\begin{equation}
Y_{s}^{t,\xi,u}=\mathbb{\tilde{E}}_{s}\left[  \Phi(X_{T}^{t,\xi,u})+\int
_{s}^{T}f(r,X_{r}^{t,\xi,u},Y_{r}^{t,\xi,u},u_{r})dr+\int_{s}^{T}%
g_{ij}(r,X_{r}^{t,\xi,u},Y_{r}^{t,\xi,u},u_{r})d\langle B^{i},B^{j}\rangle
_{r}\right]  , \label{e3}%
\end{equation}
where $s\in \lbrack t,T]$, $\langle B\rangle=(\langle B^{i},B^{j}%
\rangle)_{i,j=1}^{d}$ is the quadratic variation of $B$.

Suppose that $b$, $h_{ij}:[0,T]\times \mathbb{R}^{n}\times U\rightarrow
\mathbb{R}^{n}$, $\sigma:[0,T]\times \mathbb{R}^{n}\times U\rightarrow
\mathbb{R}^{n\times d}$, $\Phi:\mathbb{R}^{n}\rightarrow \mathbb{R}$, $f$,
$g_{ij}:[0,T]\times \mathbb{R}^{n}\times \mathbb{R}\times U\rightarrow
\mathbb{R}$ are deterministic functions and satisfy the following conditions:

\begin{description}
\item[(H1)] There exists a constant $L>0$ such that for any $(s,x,y,u)$,
$(s,x^{\prime},y^{\prime},v)\in \lbrack0,T]\times \mathbb{R}^{n}\times
\mathbb{R}\times U$,%
\[%
\begin{array}
[c]{l}%
|b(s,x,u)-b(s,x^{\prime},v)|+|h_{ij}(s,x,u)-h_{ij}(s,x^{\prime},v)|+|\sigma
(s,x,u)-\sigma(s,x^{\prime},v)|\\
\  \  \leq L(|x-x^{\prime}|+|u-v|),\\
|\Phi(x)-\Phi(x^{\prime})|\leq L|x-x^{\prime}|,\\
|f(s,x,y,u)-f(s,x^{\prime},y^{\prime},v)|+|g_{ij}(s,x,y,u)-g_{ij}(s,x^{\prime
},y^{\prime},v)|\\
\  \  \leq L(|x-x^{\prime}|+|y-y^{\prime}|+|u-v|);
\end{array}
\]

\item[(H2)] $h_{ij}=h_{ji}$ and $g_{ij}=g_{ji}$; $b,$ $h_{ij},$ $\sigma,$ $f,$
$g_{ij}$ are continuous in $s$.
\end{description}

We have the following theorems.

\begin{theorem}
(\cite{P2019}) Let Assumptions (H1) and (H2) hold. Then, for each $\xi \in
L_{G}^{2}(\Omega_{t};\mathbb{R}^{n})$ and $u\in \mathcal{U}[t,T]$, there exists
a unique solution $(X,Y)\in M_{G}^{2}(t,T;\mathbb{R}^{n+1})$ for the
forward-backward SDE (\ref{e2}) and (\ref{e3}).
\end{theorem}

\begin{theorem}
\label{th-ex}(\cite{HJ1, P2019}) Let Assumptions (H1) and (H2) hold, and let
$\xi,$ $\xi^{\prime}\in L_{G}^{p}(\Omega_{t};\mathbb{R}^{n})$ with $p\geq2$
and $u,$ $v\in \mathcal{U}[t,T]$. Then, for each $\delta \in \lbrack0,T-t]$, we
have%
\[%
\begin{array}
[c]{l}%
\mathbb{\hat{E}}_{t}[|X_{t+\delta}^{t,\xi,u}-X_{t+\delta}^{t,\xi^{\prime}%
,v}|^{2}]\leq C(|\xi-\xi^{\prime}|^{2}+\mathbb{\hat{E}}_{t}[\int_{t}%
^{t+\delta}|u_{s}-v_{s}|^{2}ds]),\\
\mathbb{\hat{E}}_{t}[|X_{t+\delta}^{t,\xi,u}|^{p}]\leq C(1+|\xi|^{p}),\\
\mathbb{\hat{E}}_{t}\left[  \underset{s\in \lbrack t,t+\delta]}{\sup}%
|X_{s}^{t,\xi,u}-\xi|^{p}\right]  \leq C(1+|\xi|^{p})\delta^{p/2},
\end{array}
\]
where $C$ depends on $T$, $G$, $p$ and $L$.
\end{theorem}

Our stochastic optimal control problem is to find $u\in \mathcal{U}[t,T]$ which
minimizes the objective function $Y_{t}^{t,x,u}$ for each given $x\in
\mathbb{R}^{n}$. For this purpose, we need the following definition of
\ essential infimum of $\{Y_{t}^{t,x,u}:u\in \mathcal{U}[t,T]\}$.

\begin{definition}
(\cite{HJ1}) The essential infimum of $\{Y_{t}^{t,x,u}:u\in \mathcal{U}%
[t,T]\}$, denoted by $\underset{u\in \mathcal{U}[t,T]}{ess\inf}Y_{t}^{t,x,u}$,
is a random variable $\zeta \in L_{G}^{2}(\Omega_{t})$ satisfying:

\begin{description}
\item[(i)] for any $u\in \mathcal{U}[t,T]$, $\zeta \leq Y_{t}^{t,x,u}$ q.s.;

\item[(ii)] if $\eta$ is a random variable satisfying $\eta \leq Y_{t}^{t,x,u}$
q.s. for any $u\in \mathcal{U}[t,T]$, then $\zeta \geq \eta$ q.s.
\end{description}
\end{definition}

For each $(t,x)\in \lbrack0,T]\times \mathbb{R}^{n}$, we define the value
function%
\begin{equation}
V(t,x):=\underset{u\in \mathcal{U}[t,T]}{ess\inf}Y_{t}^{t,x,u}. \label{e4}%
\end{equation}

In the following we will prove that $V(\cdot,\cdot)$ is deterministic and
$V(t,\xi)=\underset{u\in \mathcal{U}[t,T]}{ess\inf}Y_{t}^{t,\xi,u}$ for each
$\xi \in L_{G}^{2}(\Omega_{t};\mathbb{R}^{n})$. Furthermore, we will obtain the
dynamic programming principle and the related fully nonlinear HJB equation.

\section{Dynamic programming principle}

In the following, the constant $C$ will change from line to line in our proof.
We use the following notations: for each given $0\leq t\leq s\leq T$,%

\[%
\begin{array}
[c]{l}%
Lip(\Omega_{s}^{t}):=\{ \varphi(B_{t_{1}}-B_{t},\ldots,B_{t_{N}}-B_{t}%
):N\geq1,t_{1},\ldots,t_{N}\in \lbrack t,s],\varphi \in C_{b.Lip}(\mathbb{R}%
^{d\times N})\};\\
L_{G}^{2}(\Omega_{s}^{t}):=\{ \text{the completion of }Lip(\Omega_{s}%
^{t})\text{ under the norm }||\cdot||_{L_{G}^{2}}\};\\
M_{G}^{0,t}(t,T):=\{ \eta_{s}=\sum_{k=0}^{N-1}\xi_{k}I_{[t_{k},t_{k+1}%
)}(s):t=t_{0}<\cdots<t_{N}=T,\xi_{k}\in Lip(\Omega_{t_{k}}^{t})\};\\
M_{G}^{2,t}(t,T):=\{ \text{the completion of }M_{G}^{0,t}(t,T)\text{ under the
norm }||\cdot||_{M_{G}^{2}}\};\\
\mathcal{U}^{t}[t,T]:=\{u:u\in M_{G}^{2,t}(t,T;\mathbb{R}^{m})\text{ with
values in }U\};\\
\mathbb{U}[t,T]:=\{u=\sum_{k=1}^{N}I_{A_{k}}u^{k}:N\geq1,u^{k}\in
\mathcal{U}^{t}[t,T],I_{A_{k}}\in L_{G}^{2}(\Omega_{t}),(A_{k})_{k=1}%
^{N}\text{ is a partition of }\Omega \}.
\end{array}
\]

In order to prove that $V(\cdot,\cdot)$ is deterministic, we need the
following two lemmas. The first lemma can be found in \cite{HJ1}.

\begin{lemma}
(\cite{HJ1}) Let $u\in \mathcal{U}[t,T]$ be given. Then there exists a sequence
$(u^{k})_{k\geq1}$ in $\mathbb{U}[t,T]$ such that%
\[
\lim_{k\rightarrow \infty}\mathbb{\hat{E}}\left[  \int_{t}^{T}|u_{s}-u_{s}%
^{k}|^{2}ds\right]  =0.
\]

\end{lemma}

\begin{lemma}
Let Assumptions (H1) and (H2) hold, and let $\xi \in L_{G}^{2}(\Omega
_{t};\mathbb{R}^{n})$, $u\in \mathcal{U}[t,T]$ and $v=\sum_{k=1}^{N}I_{A_{k}%
}v^{k}\in \mathbb{U}[t,T]$. Then there exists a constant $C$ depending on $T$,
$G$ and $L$ such that%
\[
\mathbb{\hat{E}}\left[  \left \vert Y_{t}^{t,\xi,u}-\sum_{k=1}^{N}I_{A_{k}%
}Y_{t}^{t,\xi,v^{k}}\right \vert ^{2}\right]  \leq C\mathbb{\hat{E}}\left[
\int_{t}^{T}|u_{s}-v_{s}|^{2}ds\right]  .
\]

\end{lemma}

\begin{proof}
Similar to the proof of Lemma 15 in \cite{HJ1}, we can get%
\[
X_{s}^{t,\xi,v}=\sum_{k=1}^{N}I_{A_{k}}X_{s}^{t,\xi,v^{k}}\text{ and }%
Y_{s}^{t,\xi,v}=\sum_{k=1}^{N}I_{A_{k}}Y_{s}^{t,\xi,v^{k}}\text{ for }%
s\in \lbrack t,T].
\]
Since $\tilde{G}$-expectation $\mathbb{\tilde{E}}$ is dominated by
$G$-expectation $\mathbb{\hat{E}}$, by (\ref{e3}), we obtain
\[
|Y_{s}^{t,\xi,u}-Y_{s}^{t,\xi,v}|\leq C\mathbb{\hat{E}}_{s}\left[
|X_{T}^{t,\xi,u}-X_{T}^{t,\xi,v}|+\int_{s}^{T}(|Y_{r}^{t,\xi,u}-Y_{r}%
^{t,\xi,v}|+|X_{r}^{t,\xi,u}-X_{r}^{t,\xi,v}|+|u_{r}-v_{r}|)dr\right]  ,
\]
where $s\in \lbrack t,T]$ and $C$ depends on $G$ and $L$. By the H\"{o}lder
inequality, we get
\[
|Y_{s}^{t,\xi,u}-Y_{s}^{t,\xi,v}|^{2}\leq C\mathbb{\hat{E}}_{s}\left[
|X_{T}^{t,\xi,u}-X_{T}^{t,\xi,v}|^{2}+\int_{s}^{T}(|Y_{r}^{t,\xi,u}%
-Y_{r}^{t,\xi,v}|^{2}+|X_{r}^{t,\xi,u}-X_{r}^{t,\xi,v}|^{2}+|u_{r}-v_{r}%
|^{2})dr\right]  ,
\]
where $s\in \lbrack t,T]$ and $C$ depends on $T$, $G$ and $L$. By the Gronwall
inequality under $\mathbb{\hat{E}}$ (see Theorem 3.10 in \cite{HJPS}), we
deduce%
\begin{equation}
|Y_{t}^{t,\xi,u}-Y_{t}^{t,\xi,v}|^{2}\leq C\mathbb{\hat{E}}_{t}\left[
|X_{T}^{t,\xi,u}-X_{T}^{t,\xi,v}|^{2}+\int_{t}^{T}(|X_{r}^{t,\xi,u}%
-X_{r}^{t,\xi,v}|^{2}+|u_{r}-v_{r}|^{2})dr\right]  , \label{e5}%
\end{equation}
where $C$ depends on $T$, $G$ and $L$. By Theorem \ref{th-ex}, we have
\begin{equation}
\mathbb{\hat{E}}_{t}\left[  |X_{s}^{t,\xi,u}-X_{s}^{t,\xi,v}|^{2}\right]  \leq
C\mathbb{\hat{E}}_{t}\left[  \int_{t}^{T}|u_{r}-v_{r}|^{2}dr\right]  ,
\label{e6}%
\end{equation}
where $s\in \lbrack t,T]$ and $C$ depends on $T$, $G$ and $L$. Thus we obtain
the desired result by (\ref{e5}) and (\ref{e6}).
\end{proof}

\begin{theorem}
\label{th-vtx}Let Assumptions (H1) and (H2) hold. Then the value function
$V(t,x)$ exists and%
\[
V(t,x)=\inf_{u\in \mathcal{U}^{t}[t,T]}Y_{t}^{t,x,u}.
\]

\end{theorem}

\begin{proof}
The proof is the same as Theorem 17 in \cite{HJ1}. We omit it.
\end{proof}

Now we study the properties of $V(\cdot,\cdot)$.

\begin{proposition}
\label{pro-v}Let Assumptions (H1) and (H2) hold. Then there exists a constant
$C$ depending on $T$, $G$ and $L$ such that, for any $t\in \lbrack0,T]$, $x,$
$y\in \mathbb{R}^{n},$%
\[
|V(t,x)-V(t,y)|\leq C|x-y|\text{ and }|V(t,x)|\leq C(1+|x|).
\]

\end{proposition}

\begin{proof}
Similar to the proof of inequality (\ref{e5}), we can obtain that, for any
$u\in \mathcal{U}^{t}[t,T]$,%
\begin{equation}
|Y_{t}^{t,x,u}-Y_{t}^{t,y,u}|^{2}\leq C\mathbb{\hat{E}}_{t}\left[
|X_{T}^{t,x,u}-X_{T}^{t,y,u}|^{2}+\int_{t}^{T}|X_{r}^{t,x,u}-X_{r}%
^{t,y,u}|^{2}dr\right]  , \label{e7}%
\end{equation}
where $C$ depends on $T$, $G$ and $L$. By Theorem \ref{th-ex}, we have%
\begin{equation}
\mathbb{\hat{E}}_{t}\left[  |X_{s}^{t,x,u}-X_{s}^{t,y,u}|^{2}\right]  \leq
C|x-y|^{2}, \label{e8}%
\end{equation}
where $s\in \lbrack t,T]$ and $C$ depends on $T$, $G$ and $L$. Thus we get
$|V(t,x)-V(t,y)|\leq C|x-y|$ by (\ref{e7}) and (\ref{e8}). Similarly, we can
obtain $|V(t,x)|\leq C(1+|x|)$.
\end{proof}

\begin{theorem}
\label{th-vxi}Let Assumptions (H1) and (H2) hold. Then, for any $\xi \in
L_{G}^{2}(\Omega_{t};\mathbb{R}^{n})$, we have%
\[
V(t,\xi)=\underset{u\in \mathcal{U}[t,T]}{ess\inf}Y_{t}^{t,\xi,u}.
\]

\end{theorem}

\begin{proof}
For each given $u\in \mathcal{U}[t,T]$, we first prove that $V(t,\xi)\leq
Y_{t}^{t,\xi,u}$ q.s.

For each $\varepsilon>0$, we can find a $\xi_{\varepsilon}=\sum_{k=1}^{\infty
}x_{k}I_{A_{k}}$ such that $|\xi-\xi_{\varepsilon}|\leq \varepsilon$, where
$x_{k}\in \mathbb{R}^{n}$ and $\{A_{k}\}_{k=1}^{\infty}$ is a $\mathcal{B}%
(\Omega_{t})$-partition of $\Omega$. By Proposition \ref{pro-v}, we get%
\begin{equation}
\left \vert V(t,\xi)-\sum_{k=1}^{\infty}V(t,x_{k})I_{A_{k}}\right \vert
=\left \vert V(t,\xi)-V(t,\xi_{\varepsilon})\right \vert \leq C\varepsilon.
\label{e9}%
\end{equation}
Similar to the proof of inequalities (\ref{e7}) and (\ref{e8}), we can get%
\[
|Y_{t}^{t,\xi,u}-Y_{t}^{t,x_{k},u}|\leq C|\xi-x_{k}|,\text{ }k\geq1,
\]
where $C$ depends on $T$, $G$ and $L$. Then, we obtain%
\begin{equation}
\left \vert Y_{t}^{t,\xi,u}-\sum_{k=1}^{\infty}Y_{t}^{t,x_{k},u}I_{A_{k}%
}\right \vert =\sum_{k=1}^{\infty}|Y_{t}^{t,\xi,u}-Y_{t}^{t,x_{k},u}|I_{A_{k}%
}\leq C|\xi-\xi_{\varepsilon}|\leq C\varepsilon. \label{e10}%
\end{equation}
By (\ref{e4}), we have%
\begin{equation}
\sum_{k=1}^{\infty}V(t,x_{k})I_{A_{k}}\leq \sum_{k=1}^{\infty}Y_{t}^{t,x_{k}%
,u}I_{A_{k}},\text{ q.s.} \label{e11}%
\end{equation}
It follows from (\ref{e9}), (\ref{e10}) and (\ref{e11}) that%
\[
V(t,\xi)\leq Y_{t}^{t,\xi,u}+C\varepsilon,\text{ q.s.,}%
\]
where $C$ is independent of $\varepsilon$. Thus we obtain $V(t,\xi)\leq
Y_{t}^{t,\xi,u}$ q.s.

Second, if $\eta$ is a random variable satisfying $\eta \leq Y_{t}^{t,\xi,u}$
q.s. for any $u\in \mathcal{U}[t,T]$, then we prove that $V(t,\xi)\geq \eta$ q.s.

It is easy to verify that the constant $C$ in inequality (\ref{e10}) is
independent of $u$. Then we obtain%
\[
\eta \leq \sum_{k=1}^{\infty}Y_{t}^{t,x_{k},u}I_{A_{k}}+C\varepsilon,\text{
q.s., for any }u\in \mathcal{U}[t,T],
\]
where $C$ depends on $T$, $G$ and $L$. By Theorem \ref{th-vtx} and the above
inequality, we can get%
\begin{equation}
\eta \leq \sum_{k=1}^{\infty}V(t,x_{k})I_{A_{k}}+C\varepsilon,\text{ q.s.}
\label{e12}%
\end{equation}
Thus we obtain $V(t,\xi)\geq \eta$ q.s. by (\ref{e9}) and (\ref{e12}), which
implies the desired result.
\end{proof}

Finally, we study the dynamic programming principle. The following lemma is
useful in deriving the dynamic programming principle.

\begin{lemma}
\label{le-apxi}Let $\xi \in L_{G}^{2}(\Omega_{s})$ with fixed $s\in \lbrack
0,T]$. Then there exists a sequence $\xi_{k}=\sum_{i=1}^{N_{k}}x_{i}%
^{k}I_{A_{i}^{k}}$, $k\geq1$, such that%
\[
\lim_{k\rightarrow \infty}\mathbb{\hat{E}}\left[  |\xi-\xi_{k}|^{2}\right]
=0,
\]
where $x_{i}^{k}\in \mathbb{R}$, $I_{A_{i}^{k}}\in L_{G}^{2}(\Omega_{s})$,
$i\leq N_{k}$, $k\geq1$ and $(A_{i}^{k})_{i=1}^{N_{k}}$ is a $\mathcal{B}%
(\Omega_{s})$-partition of $\Omega$.
\end{lemma}

\begin{proof}
Since $L_{G}^{2}(\Omega_{s})$ is the completion of $Lip(\Omega_{s})$ under the
norm $||\cdot||_{2}$, we only need to prove the case
\[
\xi=\varphi(B_{t_{1}},B_{t_{2}}-B_{t_{1}},\ldots,B_{t_{N}}-B_{t_{N-1}}),
\]
where $N\geq1$, $0<t_{1}<\cdots<t_{N}\leq s$, $\varphi \in C_{b.Lip}%
(\mathbb{R}^{d\times N})$.

By Theorem 3.20 in \cite{HWZ}, we know that
\[
I_{\{(B_{t_{1}},B_{t_{2}}-B_{t_{1}},\ldots,B_{t_{N}}-B_{t_{N-1}})\in \lbrack
c,c^{\prime})\}}\in L_{G}^{2}(\Omega_{s})
\]
for each $c$, $c^{\prime}\in \mathbb{R}^{d\times N}$ with $c\leq c^{\prime}$.
For each $k\geq1$, we can find
\[
A_{i}^{k}=\{(B_{t_{1}},B_{t_{2}}-B_{t_{1}},\ldots,B_{t_{N}}-B_{t_{N-1}}%
)\in \lbrack c_{i,k},c_{i,k}^{\prime})\},i=1,\ldots,N_{k}-1,
\]
such that $[-ke,ke)=\cup_{i\leq N_{k}-1}[c_{i,k},c_{i,k}^{\prime})$ with
$e=[1,\ldots,1]^{T}\in \mathbb{R}^{d\times N}$, $|c_{i,k}^{\prime}-c_{i,k}|\leq
k^{-1}$ and $A_{i}^{k}\cap A_{j}^{k}=\emptyset$ for $i\neq j$. Set $A_{N_{k}%
}^{k}=\Omega \backslash \cup_{i\leq N_{k}-1}A_{i}^{k}$ and
\[
\xi_{k}=\sum_{i=1}^{N_{k}-1}\varphi(c_{i,k})I_{A_{i}^{k}}+0I_{A_{N_{k}}^{k}}.
\]
Then we obtain%
\[
|\xi-\xi_{k}|\leq \frac{L_{\varphi}}{k}+\frac{M_{\varphi}}{k}(|B_{t_{1}%
}|+|B_{t_{2}}-B_{t_{1}}|+\cdots+|B_{t_{N}}-B_{t_{N-1}}|),
\]
where $L_{\varphi}$ is the Lipschitz constant of $\varphi$ and $M_{\varphi}$
is the bound of $\varphi$. Thus%
\[
\mathbb{\hat{E}}\left[  |\xi-\xi_{k}|^{2}\right]  \leq \frac{C}{k^{2}},
\]
which yields the desired result.
\end{proof}

In order to give the dynamic programming principle, we define the following
backward semigroup $\mathbb{G}_{t,t+\delta}^{t,x,u}[\cdot]$ which was first
introduced by Peng in \cite{peng-dpp-1}.

For each given $(t,x)\in \lbrack0,T)\times \mathbb{R}^{n}$, $\delta \in
\lbrack0,T-t]$, $u\in \mathcal{U}[t,t+\delta]$ and $\eta \in L_{G}^{2}%
(\Omega_{t+\delta})$, define%
\[
\mathbb{G}_{s,t+\delta}^{t,x,u}[\eta]=\tilde{Y}_{s}^{t,x,u}\text{ for }%
s\in \lbrack t,t+\delta],
\]
where $(X_{s}^{t,x,u},\tilde{Y}_{s}^{t,x,u})_{s\in \lbrack t,t+\delta]}$ is the
solution of the following forward and backward SDEs:
\begin{equation}
\left \{
\begin{array}
[c]{rl}%
dX_{s}^{t,x,u}= & b(s,X_{s}^{t,x,u},u_{s})ds+h_{ij}(s,X_{s}^{t,x,u}%
,u_{s})d\langle B^{i},B^{j}\rangle_{s}+\sigma(s,X_{s}^{t,x,u},u_{s})dB_{s},\\
X_{t}^{t,x,u}= & x,
\end{array}
\right.  \label{newe2}%
\end{equation}
and%
\begin{equation}
\tilde{Y}_{s}^{t,x,u}=\mathbb{\tilde{E}}_{s}\left[  \eta+\int_{s}^{t+\delta
}f(r,X_{r}^{t,x,u},\tilde{Y}_{r}^{t,x,u},u_{r})dr+\int_{s}^{t+\delta}%
g_{ij}(r,X_{r}^{t,x,u},\tilde{Y}_{r}^{t,x,u},u_{r})d\langle B^{i},B^{j}%
\rangle_{r}\right]  . \label{e13}%
\end{equation}

The following lemma is the comparison theorem of backward SDE under
$\mathbb{\tilde{E}}$.

\begin{lemma}
\label{le-com}Let Assumptions (H1) and (H2) hold, and let $(t,x)\in
\lbrack0,T)\times \mathbb{R}^{n}$, $\delta \in \lbrack0,T-t]$, $u\in
\mathcal{U}[t,t+\delta]$ and $\eta_{1}$, $\eta_{2}\in L_{G}^{2}(\Omega
_{t+\delta})$ be given. If $\eta_{1}\geq \eta_{2}$ q.s., then $\mathbb{G}%
_{t,t+\delta}^{t,x,u}[\eta_{1}]\geq \mathbb{G}_{t,t+\delta}^{t,x,u}[\eta_{2}]$ q.s.
\end{lemma}

\begin{proof}
Denote $Y_{s}^{1}=\mathbb{G}_{s,t+\delta}^{t,x,u}[\eta_{1}]$, $Y_{s}%
^{2}=\mathbb{G}_{s,t+\delta}^{t,x,u}[\eta_{2}]$, $\hat{Y}_{s}=Y_{s}^{1}%
-Y_{s}^{2}$ for $s\in \lbrack t,t+\delta]$, and $\hat{\eta}=\eta_{1}-\eta_{2}$.
For each given $\varepsilon>0$, just like the proof of Theorem 3.6 in
\cite{HJPS}, we can find $(a_{s}^{\varepsilon})_{s\in \lbrack t,t+\delta]}$,
$(m_{s}^{\varepsilon})_{s\in \lbrack t,t+\delta]}$, $(c_{s}^{ij,\varepsilon
})_{s\in \lbrack t,t+\delta]}$, $(n_{s}^{ij,\varepsilon})_{s\in \lbrack
t,t+\delta]}\in M_{G}^{2}(t,t+\delta)$ such that $|a_{s}^{\varepsilon}|\leq
L$, $|c_{s}^{ij,\varepsilon}|\leq L$, $|m_{s}^{\varepsilon}|\leq2L\varepsilon
$, $|n_{s}^{ij,\varepsilon}|\leq2L\varepsilon$,%
\[
f(r,X_{r}^{t,x,u},Y_{r}^{1},u_{r})-f(r,X_{r}^{t,x,u},Y_{r}^{2},u_{r}%
)=a_{r}^{\varepsilon}\hat{Y}_{r}+m_{r}^{\varepsilon}%
\]
and%
\[
g_{ij}(r,X_{r}^{t,x,u},Y_{r}^{1},u_{r})-g_{ij}(r,X_{r}^{t,x,u},Y_{r}^{2}%
,u_{r})=c_{r}^{ij,\varepsilon}\hat{Y}_{r}+n_{r}^{ij,\varepsilon}.
\]
Then%
\begin{equation}
\hat{Y}_{s}=\mathbb{\tilde{E}}_{s}\left[  \hat{\eta}+\tilde{\eta}+\int
_{s}^{t+\delta}(a_{r}^{\varepsilon}\hat{Y}_{r}+m_{r}^{\varepsilon})dr+\int
_{s}^{t+\delta}(c_{r}^{ij,\varepsilon}\hat{Y}_{r}+n_{r}^{ij,\varepsilon
})d\langle B^{i},B^{j}\rangle_{r}\right]  -\mathbb{\tilde{E}}_{s}[\tilde{\eta
}], \label{e15}%
\end{equation}
where $s\in \lbrack t,t+\delta]$ and $\tilde{\eta}=\eta_{2}+\int_{t}^{t+\delta
}f(r,X_{r}^{t,x,u},Y_{r}^{2},u_{r})dr+\int_{t}^{t+\delta}g_{ij}(r,X_{r}%
^{t,x,u},Y_{r}^{2},u_{r})d\langle B^{i},B^{j}\rangle_{r}$.

For each given $k\geq1$, set $t_{l}^{k}=t+l\delta k^{-1}$, $l=0$,$1$,$\ldots
$,$k$. By (\ref{e15}), one can check that, for $s\in \lbrack t_{l}^{k}%
,t_{l+1}^{k}]$, $l=k-1$,$\ldots$,$0$,%
\begin{equation}
\hat{Y}_{s}=\mathbb{\tilde{E}}_{s}\left[  \hat{Y}_{t_{l+1}^{k}}+\tilde{\eta
}+\int_{s}^{t_{l+1}^{k}}(a_{r}^{\varepsilon}\hat{Y}_{r}+m_{r}^{\varepsilon
})dr+\int_{s}^{t_{l+1}^{k}}(c_{r}^{ij,\varepsilon}\hat{Y}_{r}+n_{r}%
^{ij,\varepsilon})d\langle B^{i},B^{j}\rangle_{r}\right]  -\mathbb{\tilde{E}%
}_{s}[\tilde{\eta}]. \label{e16}%
\end{equation}
Define $(\hat{Y}_{l}^{k})_{l=0}^{n}$ backwardly as follows: set $\hat{Y}%
_{k}^{k}=\hat{\eta}$, for $l=k-1$,$\ldots$,$0$,%
\begin{equation}
\hat{Y}_{l}^{k}=\mathbb{\tilde{E}}_{t_{l}^{k}}\left[  \hat{Y}_{l+1}^{k}%
+\tilde{\eta}+\int_{t_{l}^{k}}^{t_{l+1}^{k}}(a_{r}^{\varepsilon}\hat{Y}%
_{l+1}^{k}+m_{r}^{\varepsilon})dr+\int_{t_{l}^{k}}^{t_{l+1}^{k}}%
(c_{r}^{ij,\varepsilon}\hat{Y}_{l+1}^{k}+n_{r}^{ij,\varepsilon})d\langle
B^{i},B^{j}\rangle_{r}\right]  -\mathbb{\tilde{E}}_{t_{l}^{k}}[\tilde{\eta}].
\label{e17}%
\end{equation}

Note that $|\int_{s_{1}}^{s_{2}}\zeta_{r}d\langle B^{i},B^{j}\rangle_{r}%
|\leq(\mathbb{\hat{E}}\left[  |B^{i}|^{2}\right]  \mathbb{\hat{E}}\left[
|B^{j}|^{2}\right]  )^{1/2}\int_{s_{1}}^{s_{2}}|\zeta_{r}|dr$ for each $s_{1}%
$, $s_{2}\in \lbrack t,t+\delta]$ and $\zeta \in M_{G}^{1}(t,t+\delta)$, then
one can verify that%
\[
\left \vert \int_{t_{l}^{k}}^{t_{l+1}^{k}}a_{r}^{\varepsilon}dr+\int_{t_{l}%
^{k}}^{t_{l+1}^{k}}c_{r}^{ij,\varepsilon}d\langle B^{i},B^{j}\rangle
_{r}\right \vert \leq C\int_{t_{l}^{k}}^{t_{l+1}^{k}}(|a_{r}^{\varepsilon
}|+|c_{r}^{ij,\varepsilon}|)dr\leq Ck^{-1}%
\]
and%
\[
\left \vert \int_{t_{l}^{k}}^{t_{l+1}^{k}}m_{r}^{\varepsilon}dr+\int_{t_{l}%
^{k}}^{t_{l+1}^{k}}n_{r}^{ij,\varepsilon}d\langle B^{i},B^{j}\rangle
_{r}\right \vert \leq C\int_{t_{l}^{k}}^{t_{l+1}^{k}}(|m_{r}^{\varepsilon
}|+|n_{r}^{ij,\varepsilon}|)dr\leq C\varepsilon k^{-1},
\]
where $C$ is dependent of $L$ and $\delta$ and independent of $l$. For each
$k\geq k_{0}$ with $Ck_{0}^{-1}\leq2^{-1}$, we have%
\[
\hat{Y}_{k-1}^{k}\geq \mathbb{\tilde{E}}_{t_{k-1}^{k}}[\tilde{\eta
}-C\varepsilon k^{-1}]-\mathbb{\tilde{E}}_{t_{k-1}^{k}}[\tilde{\eta
}]=-C\varepsilon k^{-1}%
\]
and%
\[
\hat{Y}_{k-2}^{k}\geq \mathbb{\tilde{E}}_{t_{k-2}^{k}}[-(1+Ck^{-1})C\varepsilon
k^{-1}+\tilde{\eta}-C\varepsilon k^{-1}]-\mathbb{\tilde{E}}_{t_{k-2}^{k}%
}[\tilde{\eta}]=-[(1+Ck^{-1})+1]C\varepsilon k^{-1}.
\]
Continuing this process, we obtain
\begin{equation}
\hat{Y}_{0}^{k}\geq-C\varepsilon k^{-1}\sum_{l=0}^{k-1}(1+Ck^{-1})^{l}%
\geq-(e^{C}-1)\varepsilon. \label{e18}%
\end{equation}

For each given $\eta \in Lip(\Omega_{t+\delta})$, define $\phi(s_{1}%
,s_{2})=\mathbb{\hat{E}}[|\mathbb{\tilde{E}}_{s_{1}}[\eta]-\mathbb{\tilde{E}%
}_{s_{2}}[\eta]|]$ for $s_{1}$, $s_{2}\in \lbrack t,t+\delta]$. By the
definition of $\mathbb{\tilde{E}}_{s}[\eta]$, one can verify that $\phi$ is a
continuous function. Then we get%
\begin{equation}
\sup_{|s_{1}-s_{2}|\leq \delta k^{-1}}\mathbb{\hat{E}}[|\mathbb{\tilde{E}%
}_{s_{1}}[\eta]-\mathbb{\tilde{E}}_{s_{2}}[\eta]|]\rightarrow0\text{ as
}k\rightarrow \infty \text{.} \label{e19}%
\end{equation}
Note that%
\[
Y_{s}^{2}=\mathbb{\tilde{E}}_{s}[\tilde{\eta}]-\int_{t}^{s}f(r,X_{r}%
^{t,x,u},Y_{r}^{2},u_{r})dr-\int_{t}^{s}g_{ij}(r,X_{r}^{t,x,u},Y_{r}^{2}%
,u_{r})d\langle B^{i},B^{j}\rangle_{r},
\]
then, by (\ref{e19}) and $\tilde{\eta}\in L_{G}^{2}(\Omega_{t+\delta})$, one
can check that%
\begin{equation}
\sup_{|s_{1}-s_{2}|\leq \delta k^{-1}}\mathbb{\hat{E}}[|Y_{s_{1}}^{2}-Y_{s_{2}%
}^{2}|]\rightarrow0\text{ as }k\rightarrow \infty \text{.} \label{e20}%
\end{equation}
Similarly, the relation (\ref{e20}) still holds for $Y^{1}$. Thus we obtain%
\begin{equation}
\gamma_{k}:=\sup_{|s_{1}-s_{2}|\leq \delta k^{-1}}\mathbb{\hat{E}}[|\hat
{Y}_{s_{1}}-\hat{Y}_{s_{2}}|]\rightarrow0\text{ as }k\rightarrow \infty \text{.}
\label{e21}%
\end{equation}

Define $\Delta_{l}^{k}=\hat{Y}_{t_{l}^{k}}-\hat{Y}_{l}^{k}$ for $l=0$%
,$1$,$\ldots$,$k$. By (\ref{e16}), (\ref{e17}) and (\ref{e21}), we get%
\begin{equation}
\mathbb{\hat{E}}[|\Delta_{l}^{k}|]\leq(1+Ck^{-1})\mathbb{\hat{E}}%
[|\Delta_{l+1}^{k}|]+Ck^{-1}\gamma_{k}, \label{e22}%
\end{equation}
where $l=k-1$,$\ldots$,$0$, $\Delta_{k}^{k}=0$, $C$ depends on $L$ and
$\delta$. Similar to (\ref{e18}), we deduce%
\begin{equation}
\mathbb{\hat{E}}[|\Delta_{0}^{k}|]=\mathbb{\hat{E}}[|\hat{Y}_{t}-\hat{Y}%
_{0}^{k}|]\leq(e^{C}-1)\gamma_{k}. \label{e23}%
\end{equation}
It follows from (\ref{e18}), (\ref{e21}) and (\ref{e23}) that $\hat{Y}_{t}%
\geq-(e^{C}-1)\varepsilon$ q.s. Since $\varepsilon$ is arbitrary, we obtain
the desired result.
\end{proof}

The following theorem is the dynamic programming principle.

\begin{theorem}
\label{DPP}Let Assumptions (H1) and (H2) hold. Then, for each $(t,x)\in
\lbrack0,T)\times \mathbb{R}^{n}$, $\delta \in \lbrack0,T-t]$, we have%
\begin{equation}
V(t,x)=\underset{u\in \mathcal{U}[t,t+\delta]}{ess\inf}\mathbb{G}_{t,t+\delta
}^{t,x,u}[V(t+\delta,X_{t+\delta}^{t,x,u})]=\inf_{u\in \mathcal{U}%
^{t}[t,t+\delta]}\mathbb{G}_{t,t+\delta}^{t,x,u}[V(t+\delta,X_{t+\delta
}^{t,x,u})]. \label{e14}%
\end{equation}

\end{theorem}

\begin{proof}
By Theorem \ref{th-vtx}, we have%
\[
\underset{u\in \mathcal{U}[t,t+\delta]}{ess\inf}\mathbb{G}_{t,t+\delta}%
^{t,x,u}[V(t+\delta,X_{t+\delta}^{t,x,u})]=\inf_{u\in \mathcal{U}%
^{t}[t,t+\delta]}\mathbb{G}_{t,t+\delta}^{t,x,u}[V(t+\delta,X_{t+\delta
}^{t,x,u})].
\]
For any $u\in \mathcal{U}^{t}[t,T]$, by Theorem \ref{th-vxi}, we get%
\[
Y_{t+\delta}^{t,x,u}=Y_{t+\delta}^{t+\delta,X_{t+\delta}^{t,x,u},u}\geq
V(t+\delta,X_{t+\delta}^{t,x,u})\text{ q.s.}%
\]
Then, by Lemma \ref{le-com}, we obtain%
\[
Y_{t}^{t,x,u}=\mathbb{G}_{t,t+\delta}^{t,x,u}[Y_{t+\delta}^{t,x,u}%
]\geq \mathbb{G}_{t,t+\delta}^{t,x,u}[V(t+\delta,X_{t+\delta}^{t,x,u})],
\]
which implies%
\[
V(t,x)\geq \inf_{u\in \mathcal{U}^{t}[t,t+\delta]}\mathbb{G}_{t,t+\delta
}^{t,x,u}[V(t+\delta,X_{t+\delta}^{t,x,u})].
\]

Now we prove the converse inequality. For each given $\varepsilon>0$, there
exists a $v\in \mathcal{U}^{t}[t,t+\delta]$ such that%
\begin{equation}
\mathbb{G}_{t,t+\delta}^{t,x,v}[V(t+\delta,X_{t+\delta}^{t,x,v})]\leq
\varepsilon+\inf_{u\in \mathcal{U}^{t}[t,t+\delta]}\mathbb{G}_{t,t+\delta
}^{t,x,u}[V(t+\delta,X_{t+\delta}^{t,x,u})]. \label{e25}%
\end{equation}
Since $X_{t+\delta}^{t,x,v}\in L_{G}^{2}(\Omega_{t+\delta}^{t})$, by Lemma
\ref{le-apxi}, we can find a sequence $\xi_{k}=\sum_{l=1}^{N_{k}}x_{l}%
^{k}I_{A_{l}^{k}}$, $k\geq1$, such that%
\begin{equation}
\mathbb{\hat{E}}\left[  |X_{t+\delta}^{t,x,v}-\xi_{k}|^{2}\right]  \leq
k^{-1}, \label{e24}%
\end{equation}
where $x_{l}^{k}\in \mathbb{R}$, $I_{A_{l}^{k}}\in L_{G}^{2}(\Omega_{t+\delta
}^{t})$, $l\leq N_{k}$, $k\geq1$ and $(A_{l}^{k})_{l=1}^{N_{k}}$ is a
$\mathcal{B}(\Omega_{t+\delta}^{t})$-partition of $\Omega$. For each
$x_{l}^{k}$, we can find $v_{l}^{k}\in \mathcal{U}^{t+\delta}[t+\delta,T]$ such
that%
\begin{equation}
V(t+\delta,x_{l}^{k})\leq Y_{t+\delta}^{t+\delta,x_{l}^{k},v_{l}^{k}}\leq
V(t+\delta,x_{l}^{k})+\varepsilon. \label{e26}%
\end{equation}
Set%
\[
v^{k}(s)=\sum_{l=1}^{N_{k}}v_{l}^{k}(s)I_{A_{l}^{k}}\text{ for }s\in \lbrack
t+\delta,T],
\]
and%
\[
u^{k}(s)=v(s)I_{[t,t+\delta)}(s)+v^{k}(s)I_{[t+\delta,T]}(s)\text{ for }%
s\in \lbrack t,T],
\]
it is easy to verify that $v^{k}\in \mathcal{U}[t+\delta,T]$ and $u^{k}%
\in \mathcal{U}^{t}[t,T]$. Thus we get%
\begin{equation}
V(t,x)\leq Y_{t}^{t,x,u^{k}}=\mathbb{G}_{t,t+\delta}^{t,x,v}[Y_{t+\delta
}^{t,x,u^{k}}]=\mathbb{G}_{t,t+\delta}^{t,x,v}[Y_{t+\delta}^{t+\delta
,X_{t+\delta}^{t,x,v},v^{k}}]. \label{e27}%
\end{equation}
Similarly to the proof of inequality (\ref{e5}), we obtain that%
\begin{equation}
\left \vert \mathbb{G}_{t,t+\delta}^{t,x,v}[Y_{t+\delta}^{t+\delta,X_{t+\delta
}^{t,x,v},v^{k}}]-\mathbb{G}_{t,t+\delta}^{t,x,v}[V(t+\delta,X_{t+\delta
}^{t,x,v})]\right \vert ^{2}\leq C\mathbb{\hat{E}}\left[  \left \vert
Y_{t+\delta}^{t+\delta,X_{t+\delta}^{t,x,v},v^{k}}-V(t+\delta,X_{t+\delta
}^{t,x,v})\right \vert ^{2}\right]  \label{e28}%
\end{equation}
and%
\begin{equation}
\mathbb{\hat{E}}\left[  \left \vert Y_{t+\delta}^{t+\delta,X_{t+\delta}%
^{t,x,v},v^{k}}-Y_{t+\delta}^{t+\delta,\xi_{k},v^{k}}\right \vert ^{2}\right]
\leq C\sup_{s\in \lbrack t+\delta,T]}\mathbb{\hat{E}}\left[  \left \vert
X_{s}^{t+\delta,X_{t+\delta}^{t,x,v},v^{k}}-X_{s}^{t+\delta,\xi_{k},v^{k}%
}\right \vert ^{2}\right]  , \label{e29}%
\end{equation}
where $C$ depends on $T$, $G$ and $L$. By Theorem \ref{th-ex}, (\ref{e24}) and
(\ref{e29}), we have%
\begin{equation}
\mathbb{\hat{E}}\left[  \left \vert Y_{t+\delta}^{t+\delta,X_{t+\delta}%
^{t,x,v},v^{k}}-Y_{t+\delta}^{t+\delta,\xi_{k},v^{k}}\right \vert ^{2}\right]
\leq C\mathbb{\hat{E}}\left[  \left \vert X_{t+\delta}^{t,x,v}-\xi
_{k}\right \vert ^{2}\right]  \leq Ck^{-1}, \label{e30}%
\end{equation}
where $C$ depends on $T$, $G$ and $L$. It is easy to check that%
\begin{equation}
Y_{t+\delta}^{t+\delta,\xi_{k},v^{k}}=\sum_{l=1}^{N_{k}}Y_{t+\delta}%
^{t+\delta,x_{l}^{k},v_{l}^{k}}I_{A_{l}^{k}}. \label{e31}%
\end{equation}
It follows from (\ref{e26}) and (\ref{e31}) that%
\begin{equation}
V(t+\delta,\xi_{k})\leq Y_{t+\delta}^{t+\delta,\xi_{k},v^{k}}\leq
V(t+\delta,\xi_{k})+\varepsilon. \label{e32}%
\end{equation}
By Proposition \ref{pro-v} and (\ref{e32}), we obtain%
\begin{equation}
\mathbb{\hat{E}}\left[  \left \vert Y_{t+\delta}^{t+\delta,\xi_{k},v^{k}%
}-V(t+\delta,X_{t+\delta}^{t,x,v})\right \vert ^{2}\right]  \leq C\left(
\varepsilon^{2}+\mathbb{\hat{E}}\left[  \left \vert X_{t+\delta}^{t,x,v}%
-\xi_{k}\right \vert ^{2}\right]  \right)  \leq C(\varepsilon^{2}+k^{-1}),
\label{e33}%
\end{equation}
where $C$ depends on $T$, $G$ and $L$. By (\ref{e25}), (\ref{e27}),
(\ref{e28}), (\ref{e30}) and (\ref{e33}), we deduce that%
\[
V(t,x)\leq C(\varepsilon+\sqrt{k^{-1}})+\inf_{u\in \mathcal{U}^{t}[t,t+\delta
]}\mathbb{G}_{t,t+\delta}^{t,x,u}[V(t+\delta,X_{t+\delta}^{t,x,u})],
\]
which implies the desired result by letting $k\rightarrow \infty$ and then
$\varepsilon \downarrow0$.
\end{proof}

\begin{remark}
In the above proof, we use Lemma \ref{le-apxi} to find $v^{k}$, which can be
used to simplify the proof of the dynamic programming principle and is easier
than the implied partition method in \cite{HJ1}.
\end{remark}

Now we use the dynamic programming principle to prove the continuity of
$V(\cdot,\cdot)$ in $t$.

\begin{lemma}
\label{newle333}Let Assumptions (H1) and (H2) hold. Then the value function
$V(\cdot,\cdot)$ is $\frac{1}{2}$ H\"{o}lder continuous in $t$.
\end{lemma}

\begin{proof}
For each $(t,x)\in \lbrack0,T)\times \mathbb{R}^{n}$, $\delta \in \lbrack0,T-t]$,
by Theorem \ref{DPP}, we get%
\begin{equation}
|V(t,x)-V(t+\delta,x)|\leq \sup_{u\in \mathcal{U}^{t}[t,t+\delta]}%
|\mathbb{G}_{t,t+\delta}^{t,x,u}[V(t+\delta,X_{t+\delta}^{t,x,u}%
)]-V(t+\delta,x)|. \label{e34}%
\end{equation}
For each given $u\in \mathcal{U}^{t}[t,t+\delta]$, by the definition of the
backward semigroup, we know $\mathbb{G}_{t,t+\delta}^{t,x,u}[V(t+\delta
,X_{t+\delta}^{t,x,u})]=Y_{t}$, where $(Y_{s})_{s\in \lbrack t,t+\delta]}$ is
the solution of the following backward SDE:%
\[
Y_{s}=\mathbb{\tilde{E}}_{s}\left[  V(t+\delta,X_{t+\delta}^{t,x,u})+\int
_{s}^{t+\delta}f(r,X_{r}^{t,x,u},Y_{r},u_{r})dr+\int_{s}^{t+\delta}%
g_{ij}(r,X_{r}^{t,x,u},Y_{r},u_{r})d\langle B^{i},B^{j}\rangle_{r}\right]  .
\]
By Assumptions (H1), (H2) and Proposition \ref{pro-v}, one can verify that%
\[
|Y_{s}-V(t+\delta,x)|\leq C\mathbb{\hat{E}}_{s}\left[  |X_{t+\delta}%
^{t,x,u}-x|+\int_{s}^{t+\delta}(1+|x|+|X_{r}^{t,x,u}|+|Y_{r}-V(t+\delta
,x)|)dr\right]  ,
\]
where $C$ depends on $T$, $G$ and $L$. It follows from the Gronwall inequality
under $\mathbb{\hat{E}}$ that%
\[
|Y_{t}-V(t+\delta,x)|\leq C\mathbb{\hat{E}}_{t}\left[  |X_{t+\delta}%
^{t,x,u}-x|+\int_{t}^{t+\delta}(1+|x|+|X_{r}^{t,x,u}|)dr\right]  ,
\]
where $C$ depends on $T$, $G$ and $L$. Since $\mathbb{\hat{E}}_{t}%
[|X_{t+\delta}^{t,x,u}-x|]\leq(\mathbb{\hat{E}}_{t}[|X_{t+\delta}%
^{t,x,u}-x|^{2}])^{1/2}$ and $\mathbb{\hat{E}}_{t}[|X_{r}^{t,x,u}%
|]\leq(\mathbb{\hat{E}}_{t}[|X_{r}^{t,x,u}|^{2}])^{1/2}$, we obtain%
\[
|\mathbb{G}_{t,t+\delta}^{t,x,u}[V(t+\delta,X_{t+\delta}^{t,x,u}%
)]-V(t+\delta,x)|\leq C(1+|x|)\sqrt{\delta}%
\]
by Theorem \ref{th-ex}, where $C$ depends on $T$, $G$ and $L$. Thus we obtain
$|V(t,x)-V(t+\delta,x)|\leq C(1+|x|)\sqrt{\delta}$ by inequality (\ref{e34}).
\end{proof}

\section{The viscosity solution of HJB equation}

The following theorem establishes the relationship between the value function
$V(\cdot,\cdot)$ and the fully nonlinear second-order partial differential
equation. For the definition of the viscosity solution, the readers can refer
to \cite{CIP}.

\begin{theorem}
\label{newth222}Let Assumptions (H1) and (H2) hold. Then the value function
$V(\cdot,\cdot)$ defined in (\ref{e4}) is the unique viscosity solution of the
following second-order partial differential equation:%
\begin{equation}
\left \{
\begin{array}
[c]{l}%
\partial_{t}V(t,x)+\underset{u\in U}{\inf}H(t,x,V,\partial_{x}V,\partial
_{xx}^{2}V,u)=0,\\
V(T,x)=\Phi(x),\text{ }x\in \mathbb{R}^{n},
\end{array}
\right.  \label{e35}%
\end{equation}
where%
\[
H(t,x,v,p,A,u)=\tilde{G}(F(t,x,v,p,A,u))+\langle p,b(t,x,u)\rangle
+f(t,x,v,u),
\]%
\[
F_{ij}(t,x,v,p,A,u)=(\sigma^{T}(t,x,u)A\sigma(t,x,u))_{ij}+2\langle
p,h_{ij}(t,x,u)\rangle+2g_{ij}(t,x,v,u),
\]
$(t,x,v,p,A,u)\in \lbrack0,T]\times \mathbb{R}^{n}\times \mathbb{R}%
\times \mathbb{R}^{n}\times \mathbb{S}_{n}\times U$, $\tilde{G}$ is defined in
(\ref{newe1}).
\end{theorem}

\begin{remark}
\label{new-re111}According to Theorem C.3.5 in \cite{P2019}, for the case
that
\[
\Phi \in C_{0}(\mathbb{R}^{n})=\{ \phi \in C(\mathbb{R}^{n}):\lim
_{|x|\rightarrow \infty}\phi(x)=0\},
\]
the viscosity solution of PDE (\ref{e35}) is unique; for the case that $\Phi$
$\in C(\mathbb{R}^{n})$ satisfying $|\Phi(x)|\leq C(1+|x|^{p})$ for some
positive constants $C$ and $p$, the meaning of uniqueness is that, for each
$\Phi_{k}\in C_{0}(\mathbb{R}^{n})$ such that $\Phi_{k}$ converges uniformly
to $\Phi$ on each compact set and $|\Phi_{k}|\leq C(1+|x|^{p})$, we have
$V^{\Phi_{k}}(t,x)\rightarrow V^{\Phi}(t,x)$ for $(t,x)\in \lbrack
0,T]\times \mathbb{R}^{n}$.
\end{remark}

In order to prove this theorem, we need the following lemmas. Let $\varphi \in
C_{b}^{2,3}([0,T]\times \mathbb{R}^{n})$ be given. Here $C_{b}^{2,3}%
([0,T]\times \mathbb{R}^{n})$ denotes the set of real-valued functions that are
continuously differentiable up to the second order (resp. third order) in
$t$-variable (resp. $x$-variable) and whose derivatives are bounded. For each
given $(t,x)\in \lbrack0,T)\times \mathbb{R}^{n}$, $\delta \in \lbrack0,T-t]$ and
$u\in \mathcal{U}[t,t+\delta]$, we consider the following BSDEs%
\begin{equation}
Y_{s}^{u}=\mathbb{\tilde{E}}_{s}\left[  \varphi(t+\delta,X_{t+\delta}%
^{t,x,u})+\int_{s}^{t+\delta}f(r,X_{r}^{t,x,u},Y_{r}^{u},u_{r})dr+\int
_{s}^{t+\delta}g_{ij}(r,X_{r}^{t,x,u},Y_{r}^{u},u_{r})d\langle B^{i}%
,B^{j}\rangle_{r}\right]  , \label{e36}%
\end{equation}%
\begin{equation}
Y_{s}^{1,u}=\mathbb{\tilde{E}}_{s}\left[  \int_{s}^{t+\delta}F_{1}%
(r,X_{r}^{t,x,u},Y_{r}^{1,u},u_{r})dr+\int_{s}^{t+\delta}F_{2}^{ij}%
(r,X_{r}^{t,x,u},Y_{r}^{1,u},u_{r})d\langle B^{i},B^{j}\rangle_{r}\right]
\label{e37}%
\end{equation}
and%
\begin{equation}
Y_{s}^{2,u}=\mathbb{\tilde{E}}_{s}\left[  \int_{s}^{t+\delta}F_{1}%
(r,x,0,u_{r})dr+\int_{s}^{t+\delta}F_{2}^{ij}(r,x,0,u_{r})d\langle B^{i}%
,B^{j}\rangle_{r}\right]  , \label{e38}%
\end{equation}
where $s\in \lbrack t,t+\delta]$, $(X_{s}^{t,x,u})_{s\in \lbrack t,t+\delta]}$
is the solution of the SDE (\ref{newe2}),%
\[
F_{1}(s,x,y,u)=\partial_{t}\varphi(s,x)+\langle b(s,x,u),\partial_{x}%
\varphi(s,x)\rangle+f(s,x,y+\varphi(s,x),u),
\]%
\[
F_{2}^{ij}(s,x,y,u)=\frac{1}{2}F_{ij}(s,x,y+\varphi(s,x),\partial_{x}%
\varphi(s,x),\partial_{xx}^{2}\varphi(s,x),u).
\]

\begin{lemma}
\label{new-1}For each $u\in \mathcal{U}[t,t+\delta]$, we have%
\[
Y_{s}^{1,u}=Y_{s}^{u}-\varphi(s,X_{s}^{t,x,u})\text{ for }s\in \lbrack
t,t+\delta].
\]

\end{lemma}

\begin{proof}
Applying It\^{o}'s formula to $\varphi(r,X_{r}^{t,x,u})$ on $[s,t+\delta]$, we
obtain that $(Y_{s}^{u}-\varphi(s,X_{s}^{t,x,u}))_{s\in \lbrack t,t+\delta]}$
satisfies the backward SDE (\ref{e37}), which implies the desired result by
the uniqueness of the solution.
\end{proof}

\begin{lemma}
\label{new-2}For each $u\in \mathcal{U}^{t}[t,t+\delta]$, we have%
\[
|Y_{t}^{1,u}-Y_{t}^{2,u}|\leq C(1+|x|^{3})\delta^{3/2},
\]
where the constant $C$ is dependent on $T$, $G$, $L$ and independent of $u$.
\end{lemma}

\begin{proof}
Noting that $\varphi \in C_{b}^{2,3}([0,T]\times \mathbb{R}^{n})$ and $U$ is
compact, one can verify that%
\[
|F_{1}(r,x,0,u_{r})|\leq C(1+|x|)\text{ and }|F_{2}^{ij}(r,x,0,u_{r})|\leq
C(1+|x|^{2}),
\]
where $C$ is dependent on $L$ and independent of $u$. Thus%
\begin{equation}
|Y_{s}^{2,u}|\leq C(1+|x|^{2})\delta \text{ for }s\in \lbrack t,t+\delta],
\label{e39}%
\end{equation}
where $C$ is dependent on $G$, $L$ and independent of $u$. Set $\hat{Y}%
_{s}=Y_{s}^{1,u}-Y_{s}^{2,u}$ for $s\in \lbrack t,t+\delta]$, by (\ref{e37})
and (\ref{e38}), we get%
\[
|\hat{Y}_{s}|\leq C\mathbb{\hat{E}}_{s}\left[  \int_{s}^{t+\delta}(\hat{F}%
_{r}+|\hat{Y}_{r}|)dr\right]  ,
\]
where $C>0$ is dependent on $G$, $L$ and independent of $u$,
\[
\hat{F}_{r}=|F_{1}(r,X_{r}^{t,x,u},Y_{r}^{2,u},u_{r})-F_{1}(r,x,0,u_{r}%
)|+|F_{2}^{ij}(r,X_{r}^{t,x,u},Y_{r}^{2,u},u_{r})-F_{2}^{ij}(r,x,0,u_{r})|.
\]
Note that $Y_{t}^{1,u}\in \mathbb{R}$ and $Y_{t}^{2,u}\in \mathbb{R}$ for each
$u\in \mathcal{U}^{t}[t,t+\delta]$, then, by the Gronwall inequality under
$\mathbb{\hat{E}}$, we obtain%
\begin{equation}
|Y_{t}^{1,u}-Y_{t}^{2,u}|\leq C\mathbb{\hat{E}}\left[  \int_{t}^{t+\delta}%
\hat{F}_{r}dr\right]  , \label{e40}%
\end{equation}
where $C>0$ is dependent on $T$, $G$, $L$ and independent of $u$. One can
check that%
\begin{equation}
\hat{F}_{r}\leq C\left[  (1+|x|^{2})|X_{r}^{t,x,u}-x|+|X_{r}^{t,x,u}%
-x|^{2}+|Y_{r}^{2,u}|\right]  , \label{e41}%
\end{equation}
where $C$ is dependent on $L$ and independent of $u$. It follows from
(\ref{e39}), (\ref{e40}), (\ref{e41}) and Theorem \ref{th-ex} that%
\begin{align*}
|Y_{t}^{1,u}-Y_{t}^{2,u}|  &  \leq C\left \{  (1+|x|^{2})\delta \left(
\mathbb{\hat{E}}\left[  \sup_{r\in \lbrack t,t+\delta]}|X_{r}^{t,x,u}%
-x|^{2}\right]  \right)  ^{1/2}+\delta \mathbb{\hat{E}}\left[  \sup
_{r\in \lbrack t,t+\delta]}|X_{r}^{t,x,u}-x|^{2}\right]  +(1+|x|^{2})\delta
^{2}\right \} \\
&  \leq C(1+|x|^{3})\delta^{3/2},
\end{align*}
where $C$ is dependent on $T$, $G$, $L$ and independent of $u$.
\end{proof}

\begin{lemma}
\label{newlem111}Let $\eta=(\eta^{ij})_{i,j=1}^{d}\in M_{G}^{1}(0,T;\mathbb{S}%
_{d})$. Then, for each $s\leq T$, we have%
\[
\mathbb{\tilde{E}}_{s}\left[  \int_{s}^{T}\eta_{r}^{ij}d\langle B^{i}%
,B^{j}\rangle_{r}-\int_{s}^{T}\tilde{G}(2\eta_{r})dr\right]  =0.
\]

\end{lemma}

\begin{proof}
For each $\eta$, $\tilde{\eta}\in M_{G}^{1}(0,T;\mathbb{S}_{d})$, one can
verify that
\begin{align*}
&  \mathbb{\hat{E}}\left[  \left \vert \mathbb{\tilde{E}}_{s}\left[  \int
_{s}^{T}\eta_{r}^{ij}d\langle B^{i},B^{j}\rangle_{r}-\int_{s}^{T}\tilde
{G}(2\eta_{r})dr\right]  -\mathbb{\tilde{E}}_{s}\left[  \int_{s}^{T}%
\tilde{\eta}_{r}^{ij}d\langle B^{i},B^{j}\rangle_{r}-\int_{s}^{T}\tilde
{G}(2\tilde{\eta}_{r})dr\right]  \right \vert \right] \\
&  \leq C\mathbb{\hat{E}}\left[  \int_{s}^{T}|\eta_{r}-\tilde{\eta}%
_{r}|dr\right]  ,
\end{align*}
where $C$ only depends on $G$. Thus we only need to prove the case $\eta \in
M_{G}^{0}(0,T;\mathbb{S}_{d})$, i.e.,%
\[
\eta_{r}=\sum_{k=0}^{N-1}\eta_{t_{k}}I_{[t_{k},t_{k+1})}(r),
\]
where $s=t_{0}<\cdots<t_{N}=T$, $\eta_{t_{k}}\in Lip(\Omega_{t_{k}}%
;\mathbb{S}_{d})$. Since $\mathbb{\tilde{E}}_{s}[\cdot]=\mathbb{\tilde{E}}%
_{s}[\mathbb{\tilde{E}}_{t_{k}}[\cdot]]$, we only need to prove%
\begin{equation}
\mathbb{\tilde{E}}_{t_{k}}\left[  \eta_{t_{k}}^{ij}(\langle B^{i},B^{j}%
\rangle_{t_{k+1}}-\langle B^{i},B^{j}\rangle_{t_{k}})-\tilde{G}(2\eta_{t_{k}%
})(t_{k+1}-t_{k})\right]  =0. \label{e42}%
\end{equation}

Applying It\^{o}'s formular to $\langle \eta_{t_{k}}(B_{r}-B_{t_{k}}%
),B_{r}-B_{t_{k}}\rangle$ on $[t_{k},t_{k+1}]$, we get%
\[
\mathbb{\tilde{E}}_{t_{k}}\left[  \eta_{t_{k}}^{ij}(\langle B^{i},B^{j}%
\rangle_{t_{k+1}}-\langle B^{i},B^{j}\rangle_{t_{k}})\right]  =\mathbb{\tilde
{E}}_{t_{k}}\left[  \langle \eta_{t_{k}}(B_{t_{k+1}}-B_{t_{k}}),B_{t_{k+1}%
}-B_{t_{k}}\rangle \right]  .
\]
For each given $A\in \mathbb{S}_{d}$, define
\[
u(t,x)=\mathbb{\tilde{E}}\left[  \langle A(x+B_{t}),x+B_{t}\rangle \right]
\text{ for }(t,x)\in \lbrack0,\infty)\times \mathbb{R}^{d}.
\]
By Theorem C.3.5 in \cite{P2019}, we know that $u$ is a viscosity solution of
the following PDE%
\begin{equation}
\partial_{t}u-\tilde{G}(\partial_{xx}^{2}u)=0,\text{ }u(0,x)=\langle
Ax,x\rangle. \label{e43}%
\end{equation}
On the other hand, by the proof of Theorem 3.8.2 in \cite{P2019}, we have%
\begin{equation}
u(t,x)=\langle Ax,x\rangle+\mathbb{\tilde{E}}\left[  \langle AB_{t}%
,B_{t}\rangle \right]  =\langle Ax,x\rangle+\mathbb{\tilde{E}}\left[  \langle
AB_{1},B_{1}\rangle \right]  t. \label{e44}%
\end{equation}
By (\ref{e43}) and (\ref{e44}), we obtain $\mathbb{\tilde{E}}\left[  \langle
AB_{1},B_{1}\rangle \right]  =\tilde{G}(2A)$, which implies $\mathbb{\tilde{E}%
}\left[  \langle AB_{t},B_{t}\rangle \right]  =\tilde{G}(2A)t$. Thus we have%
\[
\mathbb{\tilde{E}}_{t_{k}}\left[  \eta_{t_{k}}^{ij}(\langle B^{i},B^{j}%
\rangle_{t_{k+1}}-\langle B^{i},B^{j}\rangle_{t_{k}})\right]  =\tilde{G}%
(2\eta_{t_{k}})(t_{k+1}-t_{k}),
\]
which implies (\ref{e42}).
\end{proof}

\begin{remark}
It is important to note that we can not derive $\mathbb{\tilde{E}}\left[
\langle AB_{1},B_{1}\rangle \right]  =\tilde{G}(2A)$ by $u(t,x)=\langle
Ax,x\rangle+\tilde{G}(2A)t$ satisfying (\ref{e43}). Because, in this case of
$u(0,x)=\langle Ax,x\rangle \not \in C_{0}(\mathbb{R}^{n})$, the meaning of
uniqueness of viscosity solution is stated as in Remark \ref{new-re111}.
\end{remark}

\begin{lemma}
\label{new-3}We have%
\[
\inf_{u\in \mathcal{U}^{t}[t,t+\delta]}Y_{t}^{2,u}=\int_{t}^{t+\delta}%
F_{0}(r,x)dr,
\]
where%
\[
F_{0}(r,x)=\inf_{v\in U}\{F_{1}(r,x,0,v)+\tilde{G}(2(F_{2}^{ij}%
(r,x,0,v))_{ij=1}^{d})\}.
\]

\end{lemma}

\begin{proof}
For each $u\in \mathcal{U}^{t}[t,t+\delta]$, by Lemma \ref{newlem111}, we get
\begin{align*}
Y_{t}^{2,u}  &  =\mathbb{\tilde{E}}_{t}\left[  \int_{t}^{t+\delta}%
F_{1}(r,x,0,u_{r})dr+\int_{t}^{t+\delta}F_{2}^{ij}(r,x,0,u_{r})d\langle
B^{i},B^{j}\rangle_{r}\right] \\
&  \geq \mathbb{\tilde{E}}_{t}\left[  \int_{t}^{t+\delta}F_{0}(r,x)dr+\int
_{t}^{t+\delta}F_{2}^{ij}(r,x,0,u_{r})d\langle B^{i},B^{j}\rangle_{r}-\int
_{t}^{t+\delta}\tilde{G}(2(F_{2}^{ij}(r,x,0,u_{r}))_{ij=1}^{d})dr\right] \\
&  =\int_{t}^{t+\delta}F_{0}(r,x)dr.
\end{align*}
Hence, $\inf_{u\in \mathcal{U}^{t}[t,t+\delta]}Y_{t}^{2,u}\geq \int
_{t}^{t+\delta}F_{0}(r,x)dr$. On the other hand, we can choose a deterministic
control $u^{\ast}\in \mathcal{U}^{t}[t,t+\delta]$ such that%
\[
\int_{t}^{t+\delta}[F_{1}(r,x,0,u_{r}^{\ast})+\tilde{G}(2(F_{2}^{ij}%
(r,x,0,u_{r}^{\ast}))_{ij=1}^{d})]dr=\int_{t}^{t+\delta}F_{0}(r,x)dr.
\]
Then we obtain $Y_{t}^{2,u^{\ast}}=\int_{t}^{t+\delta}F_{0}(r,x)dr$ by Lemma
\ref{newlem111}, which implies $\inf_{u\in \mathcal{U}^{t}[t,t+\delta]}%
Y_{t}^{2,u}\leq \int_{t}^{t+\delta}F_{0}(r,x)dr$. Thus we obtain the desired result.
\end{proof}

\noindent \textbf{Proof of Theorem \ref{newth222}.} By Proposition \ref{pro-v}
and Lemma \ref{newle333}, we know that $V(\cdot,\cdot)$ is continuous on
$[0,T]\times \mathbb{R}^{n}$. Now, we first prove that $V(\cdot,\cdot)$ is the
viscosity subsolution of (\ref{e35}).

For each given $(t,x)\in \lbrack0,T)\times \mathbb{R}^{n}$, suppose $\varphi \in
C_{b}^{2,3}([0,T]\times \mathbb{R}^{n})$ such that $\varphi(t,x)=V(t,x)$ and
$\varphi \geq V$ on $[0,T]\times \mathbb{R}^{n}$. For each $\delta \in
\lbrack0,T-t]$, by Theorem \ref{DPP}, we get%
\[
V(t,x)=\inf_{u\in \mathcal{U}^{t}[t,t+\delta]}\mathbb{G}_{t,t+\delta}%
^{t,x,u}[V(t+\delta,X_{t+\delta}^{t,x,u})].
\]
Since $\varphi(t+\delta,X_{t+\delta}^{t,x,u})\geq V(t+\delta,X_{t+\delta
}^{t,x,u})$, by Lemma \ref{le-com}, we obtain $\mathbb{G}_{t,t+\delta}%
^{t,x,u}[V(t+\delta,X_{t+\delta}^{t,x,u})]\leq Y_{t}^{u}$. It follows from
$\varphi(t,x)=V(t,x)$, Lemmas \ref{new-1} and \ref{new-2} that%
\begin{align*}
\inf_{u\in \mathcal{U}^{t}[t,t+\delta]}Y_{t}^{2,u}  &  \geq \inf_{u\in
\mathcal{U}^{t}[t,t+\delta]}Y_{t}^{1,u}-C(1+|x|^{3})\delta^{3/2}\\
&  =\inf_{u\in \mathcal{U}^{t}[t,t+\delta]}(Y_{t}^{u}-\varphi(t,x))-C(1+|x|^{3}%
)\delta^{3/2}\\
&  \geq-C(1+|x|^{3})\delta^{3/2},
\end{align*}
where $C$ is dependent on $T$, $G$, $L$. By Lemma \ref{new-3}, we get%
\[
\delta^{-1}\int_{t}^{t+\delta}F_{0}(r,x)dr\geq-C(1+|x|^{3})\delta^{1/2}.
\]
One can verify that $F_{0}(\cdot,x)$ is continuous in $r$. Hence we obtain
$F_{0}(t,x)\geq0$ by letting $\delta \downarrow0$, which implies that
$V(\cdot,\cdot)$ is the viscosity subsolution of (\ref{e35}). By the same
method, we can prove that $V(\cdot,\cdot)$ is the viscosity supersolution of
(\ref{e35}). Thus $V(\cdot,\cdot)$ is the viscosity solution of (\ref{e35}).

For the uniqueness of viscosity solution, we only need to prove the case
$\Phi \in C_{0}(\mathbb{R}^{n})$ according to Remark \ref{new-re111}.
Hoverever, by the proof of Theorem C.2.9 with $l=0$ in \cite{P2019}, we see
that in order to get the uniqueness we just need to know that $\inf_{u\in
U}H(t,x,v,p,A,u)$ satisfies assumption (G$^{\prime}$). For each $t\in
\lbrack0,T)$, $x$, $y\in \mathbb{R}^{n}$, $v\in \mathbb{R}$, $\alpha>0$, $A$,
$B\in \mathbb{S}_{n}$ such that%
\[
\left(
\begin{array}
[c]{cc}%
A & 0\\
0 & B
\end{array}
\right)  \leq3\alpha \left(
\begin{array}
[c]{cc}%
I_{n} & -I_{n}\\
-I_{n} & I_{n}%
\end{array}
\right)  ,
\]
we have%
\begin{align*}
&  \underset{u\in U}{\inf}H(t,x,v,\alpha(x-y),A,u)-\underset{u\in U}{\inf
}H(t,y,v,\alpha(x-y),-B,u)\\
&  \leq \underset{u\in U}{\sup}\left[  H(t,x,v,\alpha(x-y),A,u)-H(t,y,v,\alpha
(x-y),-B,u)\right] \\
&  \leq \underset{u\in U}{\sup}G(F(t,x,v,\alpha(x-y),A,u)-F(t,y,v,\alpha
(x-y),-B,u))+L(|x-y|+\alpha|x-y|^{2})\\
&  \leq \underset{u\in U}{\sup}G(\sigma^{T}(t,x,u)A\sigma(t,x,u)+\sigma
^{T}(t,y,u)B\sigma(t,y,u))+C(|x-y|+\alpha|x-y|^{2})\\
&  \leq \underset{u\in U}{\sup}G(3\alpha(\sigma(t,x,u)-\sigma(t,y,u))^{T}%
(\sigma(t,x,u)-\sigma(t,y,u)))+C(|x-y|+\alpha|x-y|^{2})\\
&  \leq C(|x-y|+\alpha|x-y|^{2}),
\end{align*}
where $C$ depends on $L$ and $G$. Thus $\inf_{u\in U}H(t,x,v,p,A,u)$ satisfies
assumption (G$^{\prime}$), which implies that $V(\cdot,\cdot)$ is the unique
viscosity solution of (\ref{e35}). $\Box$


\bigskip

\begin{thebibliography}{99}                                                                                               %


\bibitem {BH}R. Buckdahn, Y. Hu, Probabilistic interpretation of a coupled
system of Hamilton-Jacobi-Bellman equations, J. Evol. Equ., 10 (2010), 529-549.

\bibitem {BL}R. Buckdahn, J. Li, \emph{Stochastic differential games and
viscosity solutions for Hamilton--Jacobi--Bellman--Isaacs equations}, SIAM J.
Control Optim., 47 (2008), 444--475.

\bibitem {CIP}M.G. Crandall, H. Ishii, P.L. Lions, \emph{User's guide to
viscosity solutions of second order partial differential equations, }Bull.
Amer. Math. Soc., 27 (1992), 1-67.

\bibitem {DHP11}L. Denis, M. Hu, S. Peng, \emph{Function spaces and capacity
related to a sublinear expectation: application to $G$-Brownian motion paths,}
Potential Anal., 34 (2011), 139-161.

\bibitem {DenisMartini2006}L. Denis, C. Martini, \emph{A theoretical framework
for the pricing of contingent claims in the presence of model uncertainty},
Ann. Appl. Probab., 16 (2006), 827-852.

\bibitem {DK}L. Denis, K. Kervarec, Optimal investment under model uncertainty
in non-dominated models, SIAM J. Control Optim., 51 (2013), 1803-1822.

\bibitem {EPQ}N. El Karoui, S. Peng, M. C. Quenez, \emph{Backward stochastic
differential equations in finance}, Math. Finance, 7 (1997), 1-71.

\bibitem {EJ-1}L. Epstein, S. Ji, \emph{Ambiguous volatility, possibility and
utility in continuous time,} J. Math. Econom., 50 (2014), 269-282.

\bibitem {EJ-2}L. Epstein, S. Ji, \emph{Ambiguous volatility and asset pricing
in continuous time}, Rev. Financ. Stud., 26 (2013), 1740-1786.

\bibitem {HJ1}M. Hu, S. Ji, \emph{Dynamic programming principle for stochastic
recursive optimal control problem driven by a G-Brownian motion}, Stochastic
Process. Appl., 127 (2017), 107-134.

\bibitem {HJPS1}M. Hu, S. Ji, S. Peng, Y. Song, \emph{Backward stochastic
differential equations driven by G-Brownian motion}, Stochastic Process.
Appl., 124 (2014), 759-784.

\bibitem {HJPS}M. Hu, S. Ji, S. Peng, Y. Song, \emph{Comparison theorem,
Feynman-Kac formula and Girsanov transformation for BSDEs driven by G-Brownian
motion}, Stochastic Process. Appl., 124 (2014), 1170-1195.

\bibitem {HJX}M. Hu, S. Ji, X. Xue, The existence and uniqueness of viscosity
solution to a kind of Hamilton-Jacobi-Bellman equation, SIAM J. Control
Optim., 57 (2019), 3911-3938.

\bibitem {HP09}M. Hu, S. Peng, \emph{On representation theorem of
G-expectations and paths of $G$-Brownian motion}, Acta Math. Appl. Sin. Engl.
Ser., 25 (2009), 539-546.

\bibitem {HWZ}M. Hu, F. Wang, G. Zheng, Quasi-continuous random variables and
processes under the \emph{$G$-expectation framework, }Stochastic Process.
Appl., 126 (2016), 2367-2387.

\bibitem {LW}J. Li, Q. Wei, Optimal control problems of fully coupled FBSDEs
and viscosity solutions of Hamilton-Jacobi-Bellman equations, SIAM J. Control
Optim., 52 (2014), 1622-1662.

\bibitem {MPZ}A. Matoussi, D. Possamai, C. Zhou, Robust Utility maximization
in non-dominated models with 2BSDEs, Math. Finance, 25 (2015), 258-287.

\bibitem {MPY}J. Ma, P. Protter, J. Yong, Solving forward-backward stochastic
differential equations explicitly- a four step scheme, Probab. Theory Related
Fields, 98 (1994), 339-359.

\bibitem {MY}J. Ma, J. Yong, Forward-Backward Stochastic Differential
Equations and Their Applications, Lect. Notes Math., Springer (1999).

\bibitem {Peng2005}S. Peng, \emph{Nonlinear expectations and nonlinear Markov
chains,} Chin. Ann. Math., 26B (2005), 159--184.

\bibitem {P07a}S. Peng, \emph{$G$-expectation, $G$-Brownian Motion and Related
Stochastic Calculus of It\^{o} type}, Stochastic analysis and applications,
Abel Symp., Vol. 2, Springer, Berlin, 2007, 541-567.

\bibitem {P08a}S. Peng, \emph{Multi-dimensional $G$-Brownian motion and
related stochastic calculus under $G$-expectation}, Stochastic Process. Appl.,
118 (2008), 2223-2253.

\bibitem {peng-dpp}S. Peng, \emph{A generalized dynamic programming principle
and Hamilton-Jacobi-Bellmen equation}, Stochastics Stochastics Rep., 38
(1992), 119--134.

\bibitem {peng-dpp-1}S. Peng, \emph{Backward stochastic differential
equations---stochastic optimization theory and viscosity solutions of HJB
equations}, in Topics on Stochastic Analysis, J. Yan, S. Peng, S. Fang, and L.
Wu, eds., Science Press, Beijing, 1997, 85--138 (in Chinese).

\bibitem {P2019}S. Peng, Nonlinear Expectations and Stochastic Calculus under
Uncertainty, Springer (2019).

\bibitem {PZ}T. Pham, J. Zhang, \emph{Two Person Zero-sum Game in Weak
Formulation and Path Dependent Bellman-Isaacs Equation}, SIAM J. Control
Optim., 52 (2014), 2090-2121.

\bibitem {STZ11}H. M. Soner, N. Touzi, J. Zhang, \emph{Wellposedness of Second
Order Backward SDEs,} Probab. Theory Related Fields, 153 (2012), 149-190.

\bibitem {Tang}S. Tang, Dynamic programming for general linear quadratic
optimal stochastic control with random coefficients, SIAM J. Control Optim.,
53 (2015), 1082-1106.

\bibitem {WY}Z. Wu, Z. Yu, Probabilistic interpretation for a system of
quasilinear parabolic partial differential equation combined with algebra
equations, Stochastic Process. Appl., 124 (2014), 3921-3947.

\bibitem {J.Yong}{\normalsize J. Yong, X. Y. Zhou, \emph{Stochastic controls:
Hamiltonian systems and HJB equations}, Springer (1999). }
\end{thebibliography}
\end{document}